\documentclass[12pt]{article}
\usepackage{amsxtra,amssymb,amsthm,amsmath,latexsym}

\theoremstyle{plain}

\def\oH{{\overset{\circ}{H}}}
\def\oH1{{\overset{\circ}{H}\kern-.02in{}^1}}

\def\bee{\begin{equation*}}
\def\eee{\end{equation*}}
\def\be{\begin{equation}}
\def\ee{\end{equation}}

\begin{document}

\title{Inverse obstacle scattering with
non-over-determined data
}

\author{Alexander G. Ramm\\
 Department  of Mathematics, Kansas State University, \\
 Manhattan, KS 66506, USA\\
ramm@math.ksu.edu\\
http://www.math.ksu.edu/\,$\sim $\,ramm
}

\date{}
\maketitle\thispagestyle{empty}

\begin{abstract}
\footnote{MSC: 35R30; 35J05.}
\footnote{Key words: scattering theory; obstacle scattering; uniqueness theorem; non-over-determined scattering data.
 }

It is proved that the scattering amplitude $A(\beta, \alpha_0, k_0)$, known for all $\beta\in S^2$, where $S^2$ is the
unit sphere in $\mathbb{R}^3$, and fixed $\alpha_0\in S^2$ and $k_0>0$,
 determines uniquely the surface $S$ of the obstacle $D$ and the boundary condition on $S$. The boundary condition
on $S$ is assumed to be the Dirichlet, or Neumann, or the impedance one. The uniqueness theorem for
the solution of multidimensional  inverse scattering problems with non-over-determined data was not known for many decades.
A detailed proof of such a theorem is given in this paper for inverse scattering by obstacles for the first time.  It follows
from our results that the scattering solution vanishing on the boundary $S$ of the obstacle cannot have closed surfaces of zeros in the exterior
of the obstacle different from $S$. To have a uniqueness theorem for inverse scattering problems with non-over-determined data is
of principal interest because these are the minimal scattering data that allow one to uniquely recover the scatterer.
\end{abstract}

\section{Introduction}\label{S:1}
 The uniqueness theorems for the solution of  multidimensional inverse scattering problems with non-over-determined
 scattering  data were not known since the
 origin of the inverse scattering theory, which goes, roughly speaking, to the middle of the last century.
A detailed proof of such a theorem is given in this paper for inverse scattering by obstacles for the first time.
To have a uniqueness theorem for inverse scattering problems with non-over-determined data is
of principal interest because these are the minimal scattering data that allow one to uniquely recover the scatterer.
In \cite{R584}--\cite{R603} such theorems are proved  for the first time for inverse scattering by potentials.
The result, presented in this paper was announced in \cite{R661}, where the ideas of its proof were outlined.
In this paper the arguments are given in detail, parts of the paper \cite{R661} and the ideas of its proofs are used,
two new theorems (Theorems 2 and 3) are formulated and proved, and it is pointed out that from these results
it follows that the scattering solution vanishing on the boundary $S$ of the obstacle cannot have closed surfaces of zeros
different from $S$ in the exterior of the obstacle. Theorem 1 of this paper
was the topic of my invited talk at the 2016 International Conference  on Mathematical Methods in Electromagnetic Theory
(MMET-2016), \cite{R663}.

The data is called non-over-determined if it is a function of the same number of variables as  the function to be determined from these data. In the case of the inverse scattering by an obstacle the unknown function describes the surface of this obstacle in $\mathbb{R}^3$, so
it is a function of two variables. The non-over-determined scattering data is the scattering amplitude
depending on a two-dimensional vector. The exact formulation of this inverse problem is given below.

Let us formulate the problem discussed in this paper.
Let $D\subset \mathbb{R}^3$ be a bounded domain with a connected $C^2-$smooth boundary $S$,
 $D':=\mathbb{R}^3\setminus D$ be the unbounded exterior domain and $S^2$ be the unit sphere in $\mathbb{R}^3$.
 The smoothness assumption on $S$ can be weakened.

Consider the scattering problem:
\be\label{e1} (\nabla^2+k^2)u=0 \quad in\quad D', \qquad \Gamma_j u|_{S}=0, \qquad u=e^{ik\alpha \cdot x}+v,
\ee
where the scattered field $v$ satisfies the radiation condition:
\be\label{e2}
v_r-ikv=o\Big(\frac 1 r\Big), \quad r:=|x|\to \infty.
\ee
Here $k>0$ is a constant called the wave number and $\alpha\in S^2$ is a unit vector in the direction of the propagation of the incident plane wave $e^{ik\alpha \cdot x}$. The boundary conditions are assumed to be either the Dirichlet ($\Gamma_1$), or Neumann ($\Gamma_2$),
or impedance ($\Gamma_3$) type:
\be\label{e2'}
\Gamma_1 u:=u, \quad \Gamma_2 u:=u_N, \quad \Gamma_3 u:=u_N+hu,
\ee
where  $N$ is the unit normal to $S$ pointing out of $D$,
$u_N$ is the normal derivative of $u$ on $S$,  $h=const$, Im$h\ge 0$, $h$ is the boundary impedance, and
the condition  Im$h\ge 0$ guarantees the uniqueness of the solution to the scattering problem \eqref{e1}-\eqref{e2}.

The scattering amplitude $A(\beta, \alpha,k)$ is defined by the following formula:
\be\label{e3}
v=A(\beta, \alpha,k)\frac{e^{ikr}}{r} + o\Big(\frac 1 r\Big), \quad r:=|x|\to \infty, \quad \frac {x}{r}=\beta,
\ee
where $\alpha, \beta \in S^2$, $\beta$ is the direction of the scattered wave, $\alpha$ is the direction of the incident wave.

For a bounded domain $D$ one has $ o(\frac 1 r)=O(\frac 1 {r^2})$ in formula \eqref{e3}.
The function $A(\beta, \alpha,k)$,  the scattering amplitude, can be measured experimentally. Let us call it {\em the scattering data.}
It is known (see \cite{R190}, p.25) that the solution to the scattering problem \eqref{e1}-\eqref{e2} does exist and is unique.

 {\em The inverse scattering problem (IP) consists of finding $S$ and the boundary condition on $S$ from the scattering
data.}

M.Schiffer was the first to prove in the sixties of the last century that if the boundary condition is
the Dirichlet one then the surface $S$ is uniquely determined
by the scattering data $A(\beta, \alpha_0, k)$ known for a fixed $\alpha=\alpha_0$, all $\beta\in S^2$, and all $k\in (a,b)$, $0\le a<b$.\newline
M. Schiffer  did not publish his proof. This proof can be found, for example, in \cite{R190}, p.85, and the acknowledgement of M.Schiffer's
contribution is on p.399 in \cite{R190}.

A. G. Ramm was the first to prove that the scattering data $A(\beta, \alpha, k_0)$,
known for all $\beta$ in a solid angle, all $\alpha$ in a solid angle and a fixed $k=k_0>0$ determine uniquely the boundary $S$ {\em and the
boundary condition on $S$.} This condition was assumed of one of the three types $\Gamma_j$, $j=1,2$ or $3$, (see \cite{R190}, Chapter 2, for the proof of these results).
By subindex zero fixed values of the parameters are denoted, for example, $k_0,\, \alpha_0$. By a solid angle in this paper an open subset of $S^2$ is understood.

In \cite{R190}, p.62, it is proved that for smooth bounded obstacles the scattering amplitude $A(\beta,\alpha,k)$ is an analytic
function of $\beta$ and $\alpha$ on the non-compact analytic variety \newline
$M:=\{z| z\in \mathbb{C}^3, z\cdot z=1\}$, where $z\cdot z:=\sum_{m=1}^3 z_m^2$.
The unit sphere $S^2$ is a subset of $M$. If $A(\beta,\alpha,k)$ as a function of $\beta$  is known on an open subset of $S^2$, it is uniquely extended to all of $S^2$ (and to all of $M$) by analyticity. The same is true if $A(\beta,\alpha,k)$  as a function of $\alpha$  is known on an open subset of $S^2$. By this reason one may assume that the scattering amplitude is known on all of $S^2$ if it is known in a solid angle,
that is, on  open subsets of $S^2$ as a function of $\alpha$ and $\beta$.

In papers \cite{R310} and \cite{R311} a new approach to a proof of the uniqueness theorems for inverse obstacle scattering problem (IP) was
given. This approach is used in our paper.

In paper \cite{R162} the uniqueness theorem for IP with non-over-determined data was proved for strictly convex smooth obstacles. The proof in \cite{R162}
was based on the location of resonances for a pair of such obstacles. These results are technically difficult to obtain and they hold for two strictly convex smooth obstacles with a positive distance between them.

The purpose of this paper is to prove the uniqueness theorem for IP with non-over-determined scattering data for arbitrary
$S$. For simplicity the boundary is assumed $C^2-$ smooth.
 By the boundary condition
any of the three conditions $\Gamma_j$ are understood below, but the argument is given for the Dirichlet condition for definiteness.

\vspace{3mm}

{\bf Theorem 1.}  {\em The surface $S$ and the boundary condition on $S$ are uniquely determined by the data $A(\beta)$ known
in a solid angle. }
\vspace{3mm}

{\bf Theorem 2.} {\em If  $A_1(\beta)=A_2(\beta)$ for all $\beta$ in a solid angle, then it is not possible that
$D_1\neq D_2$ and $D_1\cap D_2=\emptyset$.}

\vspace{2mm}

{\bf Theorem 3.} {\em If  $A_1(\beta)=A_2(\beta)$ for all $\beta$ in a solid angle, then it is not possible that
$D_1\neq D_2$ and $D_1\subset D_2$.}

\vspace{2mm}

{\bf Corollary.} {\em It follows from Theorems 2 and 3 that the solution to problem
\eqref{e1}-- \eqref{e2} (the scattering solution) cannot have a closed surface of zeros except
the surface $S$, the boundary of the obstacle.}

In Section 2 some auxiliary material is formulated and  Theorems 1, 2, 3 are proved. The Corollary is an immediate
consequence of these theorems. Theorem 1--3 and the Corollary are our main results.

Let us explain the logic of the proof of Theorem 1. Its proof is based on the assumption that there are
two different obstacles, $D_1$ with the surface $S_1$ and $D_2$ with the surface $S_2$,
that generate the same non-over-determined scattering data. This assumption leads to a contradiction
which proves that $S_1=S_2$. If it is proved that $S_1=S_2$, then the type of the boundary condition
(of one of the three types \eqref{e3}) can be
uniquely determined by calculating $u$  or $\frac{u_N}{u}$ on $S$.

There are three cases to consider. The first case, when $S_1$ intersects $S_2$, is considered in Theorem 1.
The second case, when $S_1$ does not intersect $S_2$ and does not lie inside $S_2$, is considered in Theorem2.
The third case, when  $S_1$ does not intersect $S_2$ and it lies inside $S_2$ is considered in Theorem 3.

 Our results show that these cases cannot occur if the non-over-determined scattering data corresponding to
$S_1$ and $S_2$ are the same. They also show that a scattering solution cannot have a closed surface of zeros except $S$.

\section{ Proofs of Theorems 1, 2 and 3}\label{S:2}

Let us formulate some lemmas which are proved by the author, except for Lemma 3, which was
known. Lemma 3 was proved first by V.Kupradze in 1934 and then by I.Vekua,  and independently
by F.Rellich, in 1943,  see a proof of Lemma 3 in the monograph \cite{R190}, p.25, and also the references there to the papers of V.Kupradze, I.Vekua, and F.Rellich).
 Another proof of Lemma 3, based on a new idea, is given in paper \cite{R136}.

Denote by $G(x,y,k)$ the Green's function corresponding to the scattering problem \eqref{e1}-\eqref{e2}. The parameter $k>0$ is
 assumed fixed in what follows. For definiteness we assume below the Dirichlet
boundary condition, but our proof is valid for the Neumann and impedance boundary conditions as well.
 If there are two surfaces $S_m$, $m=1,2$, we denote by $G_m$ the corresponding Green's functions of the Dirichlet Helmholtz operator
 in $D_m'$.

{\bf Lemma 1.} (\cite{R190}, p. 46) {\em One has:
\be\label{e4}
G(x,y,k)= g(|y|)u(x,\alpha,k) + O\Big(\frac 1 {|y|^2}\Big), \qquad |y|\to \infty, \quad \frac y{|y|}=-\alpha.
\ee
}
\vspace{2mm}

Here $g(|y|):=\frac{e^{ik|y|}}{4\pi |y|}$,  $u(x,\alpha,k)$ is the scattering solution, that is,
the solution to problem \eqref{e1}-\eqref{e2}, $ O\Big(\frac 1 {|y|^2}\Big)$ is uniform with respect to $\alpha\in S^2$,
and the notation $\gamma (r):=4\pi g(r)=\frac{e^{ik|r|}}{|r|}$ is used below.

   The solutions to equation \eqref{e1} have the unique continuation property:

  If $u$  solves  equation \eqref{e1} and vanishes on a set $\tilde{D}\subset D'$ of positive Lebesgue measure, then $u$
vanishes everywhere in $D'$.

  Formula \eqref{e4} holds
if $y$ is replaced by the vector $-\tau \alpha +\eta$,
where $\tau>0$ is a scalar and $\eta$ is an arbitrary fixed vector orthogonal to $\alpha\in S^2$, $\eta \cdot \alpha=0$.
 If  $\eta \cdot \alpha=0$ and $y=-\tau \alpha +\eta$, then $\frac {|y|}{\tau}=1+O(\frac 1 {|\tau|^2})$ as $\tau\to \infty$.
The relation  $|y|\to \infty$ is equivalent to the relation
$\tau\to \infty$, and $g(|y|)=g(\tau)(1+ O(\frac 1 {|\tau|}))$.

Denote by $D_{12}:=D_1\cup D_2$, $ D_{12}':=\mathbb{R}^3\setminus D_{12}$, $S_{12}:=\partial D_{12}$, $\tilde{S_1}:= S_{12}\setminus S_2$,
that is, $\tilde{S_1}$  does not belong to $D_2$,
$B_R':=\mathbb{R}^3\setminus B_R$, $B_R:=\{x: |x|\le R\}$. The number $R$ is sufficiently large, so that $D_{12}\subset B_R$.
Let $S^{12}$ denote the intersection of $S_1$ and $S_2$. This set may have positive two-dimensional Lebesgue measure or it may have
 two-dimensional Lebesgue measure zero.  In the first case let us denote by $L\subset S^{12}$  the line such that
 in an arbitrary small neighborhood of every point $s\in L$ there are points of $S_1$ and of $S_2$. The line $L$ has
  two-dimensional Lebesgue measure equal to zero. Denote by $S_\epsilon$ the subset of points on $S_{12}$ the distance from
  which to $L$ is less than $\epsilon$. The two-dimensional Lebesgue measure $m_\epsilon$ of  $S_\epsilon$ tends to zero as
  $\epsilon\to 0$.

A part of our proof is based on a global perturbation lemma, Lemma 2 below, which is proved in \cite{R311}, see there formula (4).
A similar lemma is proved for potential scattering in \cite{R470}, see there formula (5.1.30). For convenience of the readers a short proof of Lemma 2 is given below.

{\bf Lemma 2.}  {\em One has:
\be\label{e5}
4\pi [A_1(\beta,\alpha, k)-A_2(\beta,\alpha, k)]=\int_{S_{12}}[u_1(s,\alpha,k)u_{2N}(s,-\beta, k)-u_{1N}(s,\alpha,k)u_{2}(s,-\beta, k)]ds,
\ee
where the scattering amplitude $A_m(\beta,\alpha, k)$ corresponds to obstacle $S_m$, $m=1,2.$}
\vspace{3mm}

{\em Proof.}  Denote by $G_m(x,y,k)$ the Green's function of the Dirichlet Helmholtz operator in $D_m'$, $m=1,2$. Using Green's
formula one obtains
\be\label{e41}
G_1(x, y, k)-G_2(x,y, k)]=\int_{S_{12}}[G_1(s,x,k)G_{2N}(s, y, k)-G_{1N}(s,x,k)G_{2}(s,y, k)]ds.
\ee
Pass in \eqref{e41} to the limit $y\to \infty$, $\frac {y}{|y|}=\beta$, and use Lemma 1 to get:
\be\label{e42}
u_1(x, -\beta, k)-u_2(x, -\beta, k)]=\int_{S_{12}}[G_1(s,x,k)u_{2N}(s, -\beta, k)-G_{1N}(s,x,k)u_{2}(s,-\beta, k)]ds.
\ee
Use the formula
\be\label{e43}
u_m(x, -\beta, k)=e^{-ik\beta \cdot x}+A_m(-\alpha, -\beta, k)\frac{e^{ik|x|}}{|x|}+  O(\frac{1}{|x|^2}), \quad  |x|\to \infty, \,\, \frac{x}{|x|}=-\alpha,
\ee
pass in equation \eqref{e42} to the limit $x\to \infty$, $\frac {x}{|x|}=-\alpha$, use Lemma 1 and get
\be\label{e44}
4\pi [A_1(-\alpha,-\beta, k)-A_2(-\alpha,-\beta, k)]=\int_{S_{12}}[u_1(s,\alpha,k)u_{2N}(s,-\beta, k)-u_{1N}(s,\alpha,k)u_{2}(s,-\beta, k)]ds.
\ee
The desired relation \eqref{e5} follows from  \eqref{e44} if one recalls the known reciprocity relation
$$A(-\alpha, -\beta,k)=A(\beta, \alpha, k),$$
 which is proved, for example, in \cite{R190}, pp. 53-54.

Lemma 2 is proved. \hfill$\Box$

\vspace{3mm}
{\bf Remark 1.} In  \eqref{e41} Green's formula is used. The surface $S_{12}$ may be not smooth because
it contains the intersection $S^{12}$ of two smooth surfaces $S_1$ and $S_2$, and this intersection may be not smooth.
However, the integrand in  \eqref{e41}  is smooth up to the boundary $S_{12}$ and is uniformly bounded because  $x$ and $y$ belong to the exterior of $D_{12}$. The integral \eqref{e41} is defined as the limit of the integral over $S_{12}\setminus S_\epsilon$ as $\epsilon\to 0$ (where
$S_\epsilon$ was defined above Lemma 2). This limit does exist since $m_\epsilon$, the two-dimensional Lebesgue measure of $S_\epsilon$,
  tends to zero as $\epsilon\to 0$  while the integrand is smooth and uniformly bounded on $S_{12}$. Consequently, the integral \eqref{e41} is well defined. This argument also shows that formula \eqref{e44}  is valid for the domain $D_{12}$ if
the surfaces $S_1$ and $S_2$ are smooth and the functions $u_1$ and $u_2$ are smooth and uniformly bounded up to $S_1$ and $S_2$
respectively.

\vspace{3mm}

{\bf Lemma 3.} (\cite{R190}, p. 25) {\em If  $\lim_{r\to \infty}\int_{|x|=r}|v|^2ds=0$ and $v$ satisfies the Helmholtz
equation \eqref{e1} in $B_{R}'$, then $v=0$ in $B_R'$.}

\vspace{3mm}

The following lemma is used in our proof.

{\bf Lemma 4.} ({\em lifting lemma}) {\em If $A_1(\beta, \alpha,k)=A_2(\beta, \alpha,k)$ for all $\beta, \alpha \in S^2$, then
$G_1(x,y,k)=G_2(x,y,k)$ for all $x,y\in D_{12}'$.  If $A_1(\beta, \alpha_0, k)=A_2(\beta, \alpha_0, k)$ for all $\beta \in S^2$
and a fixed $\alpha=\alpha_0$, then $G_1(x,y_0,k)=G_2(x,y_0,k)$ for all $x\in D_{12}'$ and $y_0=-\alpha_0 \tau+\eta$, where $\tau > 0$ is a 
 number and $\eta$ is an arbitrary fixed vector orthogonal to $\alpha_0$,  $\alpha_0\cdot \eta=0$.
}

\vspace{3mm}

{\bf Proof of Lemma 4.} The function
\be\label{e45}
w:=w(x,y):=G_1(x, y, k)-G_2(x, y, k)
\ee
satisfies equation  \eqref{e1} in $D_{12}'$ as a function of $y$ and
also as a function of $x$, and $w$ satisfies the radiation condition as a function of $y$ and also as a function of $x$. By Lemma 1 one has:
\be\label{e5a}
  w=g(|y|)[u_1(x, \alpha, k)- u_2(x, \alpha, k)]+O(\frac 1 {|y|^2}), \quad |y|\to \infty, \,\, \alpha=-\frac {y}{|y|}.
 \ee
  Using formulas  \eqref{e1} and  \eqref{e3} one gets:
 \be\label{e5b}
 u_1(x, \alpha, k)- u_2(x, \alpha, k)=\gamma (|x|)[ A_1(\beta, \alpha, k)-A_2(\beta, \alpha, k)]+O(\frac 1 {|x|^2}), \quad |x|\to \infty, \,\, \beta=\frac {x}{|x|},
 \ee
because, for $m=1,2$ and $\gamma(|x|):=\frac{e^{ik|x|}}{|x|}$, one has:
\be\label{e5c}
u_m(x, \alpha, k)=e^{ik\alpha \cdot x}+A_m(\beta, \alpha, k)\gamma (|x|)+  O(\frac{1}{|x|^2}), \quad  |x|\to \infty, \,\, \beta=\frac{x}{|x|}.
\ee
If $A_1(\beta, \alpha, k)=A_2(\beta, \alpha, k)$, then equation  \eqref{e5b}  implies
\be\label{e5'}
u_1(x, \alpha, k)- u_2(x, \alpha, k)=O(\frac 1 {|x|^2}).
\ee
Since $u_1(x, \alpha, k)-u_2(x, \alpha, k)$ solves equation \eqref{e1} in $D_{12}'$ and relation \eqref{e5'} holds, it follows from Lemma 3 that $u_1(x, \alpha, k)=u_2(x, \alpha, k)$ in $B_R'$. By  the unique continuation property for the solutions to the Helmholtz equation
\eqref{e1}, one concludes that $u_1=u_2$ everywhere in $D_{12}'$. Consequently, formula  \eqref{e5a} yields
\be\label{e5d}
w(x,y)= O(\frac 1 {|y|^2}), \qquad |y|>|x|\ge R.
\ee
Since the function $w$ solves the homogeneous Helmholtz equation   \eqref{e1} in the region $|y|>|x|\ge R$, it follows
by Lemma 3 that $w=w(x,y)=0$ in this region and, by the unique continuation property, $w=0$ everywhere in $D_{12}'$.
Thus, the first part of Lemma 4 is proved.

Its second part deals with the case when $\alpha=\alpha_0$, where $\alpha_0$ is fixed. Let us prove that if
\be\label{e46}
A_1(\beta):=A_1(\beta, \alpha_0,k)=A_2(\beta, \alpha_0,k):=A_2(\beta) \qquad \forall \beta\in S^2,
\ee
then
\be\label{e47}
w(x,y_0)=0,
\ee
 where $x\in D_{12}'$ is arbitrary,  $y_0=-\tau\alpha_0 +\eta$, $\alpha_0\in S^2$ is fixed,
$\tau>0$ is a  number and $\eta$ is an arbitrary fixed vector orthogonal to $\alpha_0$,
$\eta \cdot \alpha_0=0$.

From \eqref{e46} it follows that $u_1(x,\alpha_0)=u_2(x,\alpha_0)$ for all $x\in D_{12}'$. Let us derive a contradiction
from the assumption that \eqref{e47} is not valid, or, which is equivalent, that $S_1\neq S_2$.

 The Green's formula yields $G_1(x,y_0)=g(x,y_0)-\int_{S_1}g(x,s)G_{1N}(s,y_0)ds$, where $g(x,y_0)=\frac {e^{ik|x-y_0|}}{4\pi |x-y_0|}$,
  and a similar formula holds for $G_2$  with the integration over $S_2$. Consequently,
\be\label{e47d}
G_1(x,y_0)-G_2(x,y_0)=\int_{S_2}g(s,x)G_{2N}(s,y_0)ds- \int_{S_1}g(s,x)G_{1N}(s,y_0)ds,
\ee
where  $ds$ is the surface area elements.

Let $y_0\to \infty$, $y_0/|y_0|=-\alpha_0$ and take into account that if  $u_1(x,\alpha_0)=u_2(x,\alpha_0)$ for all $x\in D_{12}'$ then $u_1(x,\alpha_0)=u_2(x,\alpha_0):=u(x,\alpha_0)$ for all $x\in D^{12}:=D_1\cap D_2$ by the unique continuation principle.
Therefore, equation \eqref{e47d} yields
\be\label{e47e}
\int_{S_2}g(s,x)u_{N}(s,\alpha_0)ds=\int_{S_1}g(s,x)u_{N}(s,\alpha_0)ds, \quad \forall x\in D^{12}.
\ee
The right side of \eqref{e47e} is an infinitely smooth function when $x$ passes the part of $S_2$ which
lies outside of $D_1$  while the normal derivative of the left side has a jump $u_{N}(s,\alpha_0)$ in such
a process. This is a contradiction unless $u_{N}(s,\alpha_0)=0$ on $S_2$. However, $u=0$ on $S_2$ and if
$ u_{N}(s,\alpha_0)=0$ on $S_2$ then, by the uniqueness of the solution to the Cauchy problem
for the Helmholtz equation, one concludes that $u=0$ in $D_2'$. This is impossible since
$\lim_{x\to \infty}|u(x,\alpha_0)|=1$. This contradiction proves that $D_1=D_2:=D$, $S_1=S_2:=S$,
and $G_1(x,y_0)=G_2(x,y_0):=G(x,y_0)$, where $G$ is the Green's function of the Dirichlet Helmholtz operator
 for the domain $D'$, and $G$ satisfies the radiation condition at infinity.

 Thus, the proof of the relation $G_1(x,y_0)=G_2(x,y_0)$
is completed  and the second part of Lemma 4 is proved.

Lemma 4 is proved. \hfill$\Box$

\vspace{2mm}

{\bf Lemma 5.} {\em One has
\be\label{e91}
\lim_{x\to t} G_{2N}(x,s,k)=\delta(s-t), \quad t\in S_2,
\ee
where $\delta(s-t)$ denotes the delta-function on $S_2$ and $x\to t$ denotes a limit along any straight line non-tangential to $S_2$.}

{\bf Proof of Lemma 5.}  Let $f\in C(S_2)$ be arbitrary. Consider the following problem: $W$ solves equation
\eqref{e1} in $D_2'$, $W$ satisfies the boundary condition $W=f$ on $S_2$, and $W$ satisfies the radiation condition.
 The unique solution to this problem is given by the Green's formula:
\be\label{e92}
W(x)=\int_{S_2} G_{2N}(x,s)f(s)ds.
\ee
 Since $\lim_{x\to t\in S_2} w(x)=f(t)$ and $f\in C(S_2)$ is arbitrary,
 the conclusion of Lemma 5 follows.

 Lemma 5 is proved. \hfill$\Box$

 \vspace {2mm}

  Let us point out the following implications:
\be\label{e5e}
G(x,y,k)\rightarrow u(x,\alpha,k)\rightarrow A(\beta, \alpha, k),
\ee
 which hold by Lemma 1 and formula  \eqref{e5c}. The first arrow means that the knowledge of $G(x,y,k)$ determines uniquely
 the scattering solution $u(x,\alpha,k)$ for all $\alpha\in S^2$, and the second arrow means that the scattering
 solution  $u(x,\alpha,k)$ determines uniquely the scattering amplitude $A(\beta, \alpha, k)$.

  The reversed implications also hold:
 \be\label{e5f}
A(\beta, \alpha, k)\rightarrow  u(x,\alpha,k)\rightarrow G(x,y, k).
\ee
 These implications follow from Lemmas 1, 3 , 4  and formula \eqref{e5c}.

 Let us explain why the knowledge of $u(x,\alpha,k)$ determines uniquely $G(x,y,k)$. If there are two $G_m$, $m=1,2,$
 to which the same $u(x,\alpha,k)$ corresponds, then $w:=G_1-G_2$ solves equation \eqref{e1} in  $D_{12}'$ and,
 by Lemma 1, $w=O(\frac 1 {|x|^2})$. Thus, by Lemma 3, $w=0$, so $G_1=G_2$ in  $D_{12}'$. This implies,
 as in the proof of Theorem 1 below, that $D_1=D_2:=D$.

Similar implications for $\alpha=\alpha_0$ fixed are formulated after the proof of Theorem 1.

\begin{proof}[Proof of Theorem 1]
If $A_1(\beta)=A_2(\beta)$ for all $\beta$ in a solid angle, then the same is true for all
$\beta\in S^2$, so one may assume that $A_1(\beta)=A_2(\beta)$ for all $\beta\in S^2$.

Let us assume that $A_1(\beta)=A_2(\beta)$ for all $\beta$ but $S_1\neq S_2$. We want to derive from this assumption
a contradiction. This contradiction will prove that the assumption  $S_1\neq S_2$ is false, so $S_1=S_2$.

If  $A_1(\beta)=A_2(\beta)$ for all $\beta\in S^2$, then Lemma 4 yields the following conclusion:
\be\label{e6}
G_1(x,y_0)=G_2(x,y_0), \qquad \forall x\in D_{12}',
\ee
where $k>0$ and $\alpha_0\in S^2$ are fixed and $y_0=-\alpha_0\tau +\eta$, $\tau>0$, $\eta\cdot \alpha_0=0$, $\eta$ is an arbitrary fixed vector
orthogonal to $\alpha_0$. This is the key point in the proof of Theorem 1. For definiteness we assume in the proof of Theorem 1 that
$S_1$ intersects $S_2$.

 If $S_1\neq S_2$ then one gets a contradiction: let $y_0$  approach a point $t\in S_2$ which does not belong to $S_1$ along the ray $-\tau \alpha_0 +\eta$.  Then, on one hand,  $G_1(x,t)=G_2(x,t)=0$ for all $x\in D_{12}'$,  and, on the other hand,
$G_1(x,t)=O(\frac 1 {|x-t|})$, so that $|G_1(x,t)|\to \infty$ as $x\to t$. This contradiction proves that $S_1=S_2$.

If $S_1=S_2:=S$ then $D_1=D_2:=D$ and  $u_1(x, \alpha_0, k)=u_2(x, \alpha_0, k):=u(x,\alpha_0,k)$ for $x\in D'$, and, consequently, the boundary condition on $S$ is uniquely determined: if $u|_{S}=0$, then one has the Dirichlet boundary condition $\Gamma_1$, otherwise
calculate $\frac{u_N}{u}$ on $S$. If this ratio vanishes, then one has the Neumann boundary condition $\Gamma_2$, otherwise one has the impedance boundary condition $\Gamma_3$, and the boundary impedance $h=-\frac{u_N}{u}$ on $S$, so the boundary condition is uniquely determined by the
non-over-determined scattering data.

Theorem 1 is proved.
\end{proof}

One may give  different proofs of Theorem 1. For example, if $S_1\neq S_2$ and $S_1$ intersects $S_2$ then, by analytic continuation,
the scattering solutions $u_m(x,\alpha_0)$, $m=1,2,$ admit analytic continuation to the exterior of the domain $D^{12}=D_1\cap D_2$.
The boundary of this domain has edges. If a point $t$ belongs to an edge, then the gradient of the solution to the homogeneous Helmholtz
equation is singular when $x\to t$. On the other hand, this $t$ belongs to a smooth boundary $S_1$ or $S_2$, so that the above gradient has to be smooth. This contradiction proves that $S_1=S_2$ in the case when  $S_1$ intersects $S_2$.

Let us formulate the implication similar to the one given before the proof of Theorem 1.
If $y = y_0 =-\tau \alpha_0 +\eta$, $\tau> 0$ is an arbitrary number, $\alpha_0$ is a fixed unit vector, and $\eta$
is an arbitrary fixed vector orthogonal to $\alpha_0$, then
\be\label{e5g}
 G(x,y_0,k)\rightarrow u(x,\alpha_0,k)\rightarrow A(\beta, \alpha_0, k),
 \ee
where  $\alpha_0$ is a free unit vector, that is, a vector whose initial point is arbitrary.

The reversed implications also hold:
 \be\label{e5h}
A(\beta, \alpha_0, k)\rightarrow  u(x,\alpha_0,k)\rightarrow G(x, y_0, k).
\ee
The first of these implications follows from Lemma 3 and the asymptotic of the scattering
solution, while the second follows from Lemmas 1 and 4.

\vspace{3mm}
We have assumed implicitly that $D_1$ and $D_2$ have a common part but none of them is a subset of the other, that is, $S-1$ intersects $S_2$.
Let us discuss the two remaining possibilities.

The first possibility is that $D_1\neq D_2$ and $D_1\cap D_2=\emptyset$.

\vspace{2mm}

{\em Proof of Theorem 2.} If $A_1(\beta)=A_2(\beta)$ in a solid angle, then  $A_1(\beta)=A_2(\beta)$ in $S^2$. This implies
that $u_1(x,\alpha_0)=u_2(x,\alpha_0)$ in $D_{12}'$. Since $u_1(x,\alpha_0)$ is defined in $D_2$ and satisfies
there the Helmholtz equation  \eqref{e1}, the unique continuation property implies that $u_2(x,\alpha_0,k)$ is defined in $D_2$
and   satisfies there the Helmholtz equation. Consequently, $u_2(x,\alpha_0,k)$ is defined in $\mathbb{R}^3$, it is
a smooth function that
satisfies  in $\mathbb{R}^3$  the Helmholtz equation, and the same is true for  $u_1(x,\alpha_0, k)$. Therefore the scattered
parts $v_1$ and $v_2$ of the scattering solutions $u_1$ and $u_2$ satisfy the Helmholtz equation \eqref{e1} in $\mathbb{R}^3$ and the radiation condition.
A function satisfying the radiation condition and the Helmholtz equation in $\mathbb{R}^3$ is equal to zero in $\mathbb{R}^3$.
Therefore, $v_1=v_2=0$ and $u_1=u_2=e^{ik\alpha_0\cdot x}$ in  $\mathbb{R}^3$. This is impossible since $u_m=0$ on
$S_m$, $m=1,2$, while $e^{ik\alpha_0\cdot x}\neq 0$ on $S_m$.

Theorem 2 is proved. \hfill$\Box$

\vspace{2mm}
The second possibility is $D_1\neq D_2$ and $D_1\subset D_2$.

\vspace{2mm}

{\em Proof of Theorem 3.}

One  proves Theorem 3 using  Lemma 4. By Lemma 4 one has
 \be\label{e32}
G_1(x,y_0)=G_2(x,y_0) \quad \forall x\in D_2',\quad y_0=-\tau \alpha_0+ \eta,\quad y_0\in D_2',\quad \eta\cdot \alpha_0=0, \quad \tau\in (0,\infty).
\ee
Note that
 \be\label{e33}
\lim_{x\to y_0}|G_1(x,y_0)|=\infty
\ee
since both $x$ and $y_0$ belong to $D_1'$ and are away from $S_1$ if $D_1\neq D_2$. On the other hand,
if $y_0\in S_2$, then $G_2(x,y_0)=0$ for all $x\in D_2$, $x\neq y_0$ and
 \be\label{e34}
\lim_{x\to y_0}|G_2(x,y_0)|=0.
\ee
This is a contradiction unless $D_1=D_2$.

Theorem 3 is proved.  \hfill$\Box$

\vspace{2mm}
{\em It follows from Theorems 1,  2 and 3 that Corollary holds: the solution to problem
\eqref{e1}-- \eqref{e2} (the scattering solution) cannot have a closed surface of zeros except
the surface $S$, the boundary of the obstacle.}

{\bf Remark 2.} Let us give
a new argument for the proof of Lemma 4.
Assume that $D_1\subset D_2$, so that geometry of
Theorem 3 is assumed.
Denote by $u$ the analytic continuation of $u_2$ into $D_2\setminus D_1$. This $u$ is equal to $u_1$ in $D_1'$.  Green's formula yields
\be\label{e81}
u(x)=\int_{S_2}g(x,s)u_N(s)ds-\int_{S_1}g(x,s)u_N(s)ds, \quad x\in D_2\setminus D_1.
\ee
Formulas \eqref{e47e} and \eqref{e81} imply 
$u=0$ in $D_2\setminus D_1$. This is a contradiction since $u$ solves the elliptic
 Helmholtz equation in $D_1'$ and if $u=0$ in  $D_2\setminus D_1$ then $u=0$ everywhere in $D_1'$, which is impossible since $u\to 1$ as $|x|\to \infty$. \hfill$\Box$

 One can argue slightly differently. In $D_2\setminus D_1$ one can derive the formulas
\be\label{e82} 
 \int_{S_2}g(x,s)u_N(s)ds=u_0(x),  \quad x\in D_2\setminus D_1.
\ee
\be\label{e83}
 \int_{S_2}g(x,s)u_N(s)ds=u_0(x)-u(x), \quad x\in D_2\setminus D_1.
 \ee 
From \eqref{e82}-\eqref{e83} it follows that $u=0$
in $D_2\setminus D_1$. This leads to a contradiction, as was explained above.

If the geometry assumed in Theorem 1 holds, then the argument is similar. The roles of $S_2$ and $S_1$ are played respectively by $S_{12}$ and 
the boundary of the intersection $D^{12}=D_1\cap D_2$.


\end{document}